\documentclass[a4paper,12pt]{article}
\usepackage{amsmath,amssymb,amsfonts,amscd,amsthm}
\usepackage{graphicx}
\usepackage{dsfont}
\usepackage{blindtext}
\newtheorem{Main}{\bg Theorem}

\usepackage[titletoc, toc, title, page]{appendix}
\usepackage{color}
\usepackage{enumerate}
\usepackage{fancyhdr}
\usepackage{parskip}		
\usepackage{lastpage}
\usepackage[utf8]{inputenc}
\usepackage[top=2cm, bottom=2cm, left=2cm, right=2cm]{geometry}
\usepackage{subcaption}
\usepackage{bbm}
\usepackage{pgf,tikz}
\usepackage{lipsum}
\usepackage{authblk}
\usetikzlibrary{arrows}
\usepackage{color}

\usepackage{xcolor}

\definecolor{bluegray}{rgb}{0.4, 0.6, 0.8}

\definecolor{dgreen}{rgb}{0.1,0.6,0.1}
\definecolor{bluegreen}{rgb}{0.1,0.5,0.2}

\def\bg{\color{bleu1}}
\def\bk{\color{black}}

\usepackage{xcolor}

\newtheorem{theorem}{\bg Theorem}                             

\newtheorem{proposition}{\bg Proposition}

\newtheorem{remark}{\bg Remark}
\newtheorem{definition}{\bg Definition}

\newcommand{\R}{\mathbb R}

\newcommand{\Z}{\mathbb Z}

\newcommand{\T}{\mathbb{T}}

\renewcommand{\a}{\alpha}

\renewcommand{\inf}{\mathop{\mbox{{inf}}}
	\renewcommand{\lim}{\mathop{\mbox{lim}}}
	\renewcommand{\max}{\mathop{\mbox{ma}x}}}



\renewcommand{\th}{\theta} 
 
\newcommand{\om}{\omega} 
\newcommand{\cO}{\mathcal O} 

\numberwithin{equation}{section}

\title{\bg \large{Topological weak
  mixing and diffusion at all times  \\ for a class of Hamiltonian systems}}
\author[1]{\small{Bassam Fayad}\thanks{Supported in part by Knut and Alice Wallenberg  foundation, grant KAW 2016.0403, and by the  ANR-15-CE40-0001.}}
\author[2]{\small{Maria Saprykina}\thanks{Supported in part by  the Swedish
    Research Council, VR 2015-04012.}}
\affil[1]{IMJ-PRG CNRS, France;  Inst. f\"or Matematik,
KTH,  Royal Institute of Technoligy, Sweden}
\affil[2]{Dept. of Mathematics, KTH Royal Institute of Technoligy, Sweden}

\date{}

\usepackage{sectsty,ulem} \definecolor{bleu1}{RGB}{0,57,128}
\subsectionfont{\color{bleu1}}
\sectionfont{\color{bleu1}}
\subsectionfont{\color{bleu1}}

\usepackage[utf8]{inputenc}
\usepackage[english]{babel}

\usepackage[colorlinks,linkcolor=bleu1,anchorcolor=green,citecolor=bleu1]{hyperref}

\usepackage{amsmath,amsthm,amssymb,mathrsfs,latexsym,amsfonts}
\usepackage{bbm}
\usepackage{amssymb}
\usepackage{amsmath}

\usepackage{color}

\theoremstyle{plain}
\newtheorem{thm}{Theorem}

\newtheorem{Lem}[thm]{\bg Lemma}

\def\Z{{\mathbb Z}}
\def\T{{\mathbb T}}
\def\R{{\mathbb R}}

\def\ZZ{{\mathbb Z}}
\def\TT{{\mathbb T}}
\def\RR{{\mathbb    R}}
\def\CC{{\mathbb C}}
\def\NN{{\mathbb N}}
\def\QQ{{\mathbb Q}}
\def\a{{\alpha}}

\newcommand{\de}{{\delta}}
\newcommand{\eps}{{\epsilon}}

\def\th{{\theta}}

\def\<{\langle}
\def\>{\rangle}

\usepackage{calc}
\usepackage{accents}

\newcommand{\pr}{{\pi}}

\newcommand{\comm}[1]{}
\newcommand{\comment}[1]{}

\usepackage{color}

\definecolor{orange}{rgb}{1,0.5,0}

\newcommand{\vertiii}[1]{{\left\vert\kern-0.25ex\left\vert\kern-0.25ex\left\vert #1 
    \right\vert\kern-0.25ex\right\vert\kern-0.25ex\right\vert}}

\begin{document}

\maketitle
\medskip 
\begin{center} {\it To our great friend and mentor Anatoly Katok}  \end{center}
\medskip 

\begin{abstract}
We present examples of nearly integrable analytic
Hamiltonian systems with several strong diffusion
properties: topological weak mixing and diffusion at all times.
These examples are obtained by AbC constructions with several frequencies.
\end{abstract}
\maketitle

\section*{Introduction}

KAM theory (after Kolmogrov, Arnol'd, and Moser) states that, 
under mild {\it non-degenera\-cy assumptions}, 
Hamiltonian systems close to integrable have their phase space almost
completely filled by invariant quasi-periodic tori. Starting from
three degrees of freedom for autonomous  Hamiltonians the existence of
such tori on an energy surface does not prevent the orbits from {\it
  circulating} between the tori inside the surface. Indeed, it was
conjectured by Arnol'd that a ``general'' 
Hamiltonian should have a dense orbit on a 
``general'' energy surface \cite{arnold}. A great amount of work has been dedicated to proving
this conjecture (giving a precise meaning to the word ``general''), 
but the picture is not yet completely clear, especially when it comes to real analytic Hamiltonians (see for example \cite{BKZ} and references therein).

In his ICM list of problems \cite{herman-ICM}, M. Herman asks: 
Can one find an example of a $C^{\infty}$-Hamiltonian  $H$ in a small
$C^k$-neighborhood, $k\geq 2$, of $H_0=\|r\|^2/2$ such that on the
energy surface $\{H=1 \}$ the Hamiltonian flow has a dense orbit?

A remarkable result in this direction is due to \cite{KZZ}: they
present an example of a Hamiltonian $H$ of the form 
$H(\th,r)=\frac{\<r,r\>}{2} +h(r,\th) \in C^{\infty}$ with a
trajectory dense in a subset of the energy surface of large measure.

Here we  present examples, with a degenerate integrable part, but with even more chaotic behavior for the perturbed system. Namely,
we show that in a class of perturbations of rotators, the generic Hamiltonian is  topologically weakly mixing on every energy surface. We also give examples of perturbed rotators for which the dynamics is diffusive at all times.

We now give the exact definitions of these properties and state our results precisely.

\begin{definition}[Topological weak mixing]
We say that the flow $\Phi^t_H$ is {\it topologically weakly mixing} if
there exists a sequence  $(t_n)_{n\in \NN}$ such that for any two open  sets $A$, $B$
on the same (arbitrarily chosen) energy surface $E_H(c)=\{(\theta,r)\mid
H(\theta,r)=c\}$,  
there exists $N=N(A,B)$ such that $\Phi^{t_n}_H(A)\cap B \neq
\emptyset$ for all $n\geq N$. 

\end{definition} 

\begin{definition}[Diffusion at all times] We say that the flow $\Phi_H^t$ exhibits {\it
diffusion at all times} if for any open set
$A \subset \TT^d\times \RR^d$ and any $R\geq 0$ there exists $T=T(A,R)$ such that
$|\pr_r(\Phi^t_H (A)) | > R$ for all $t\geq T$. Here $\pi_r$ stands for the projection onto the $r$-variables.
\end{definition}

Given $\rho>0$, denote by $C^{\om}_{\rho}$ the space
of bounded real functions on $\RR^d\times \RR^d$, that are $\ZZ^d$-periodic in the
  first $d$-vector of components, and can be extended to holomorphic
  functions on $D_{\rho}=\{(\th,r)\in (\CC^d,\CC^{d})\mid \max_j
  \{|\th_j|,|r_j|, j=1,\dots d\}\leq\rho\} $ as $\ZZ^d$-periodic with
  respect to the real part, 
with the norm $\|f\|_\rho =\sup \{|f(\th,r)|,\ (\th,r)\in D_\rho\}$.
 
Let  $\cO^{\om}_{\rho, \de}=\{f\in C^{\om}_{\rho}\mid \|f\|_\rho <\de\}$
be a neighborhood of zero in
  the space of analytic functions with the above norm. This is a Baire
  space.  
Fixing $\rho=1$ without loss of generality, we will write 
$\|f\| = \|f\|_1$, and $\cO^{\om}_{\de}$ for $\cO^{\om}_{1, \de}$.

\begin{Main}\label{th_diffusion}
There exists a dense set $Y\subset \RR^d$, $d\geq 3$, ($Y$ continuum) such
that for any $\om\in Y$, for any $\de>0$, there exists a real analytic function $h:\T^d\to \R$, such that $h \in \cO^\om_\de$, with the following property: the Hamiltonian flow
$\Phi^t_H(\theta,r)$
defined by the Hamiltonian 
$$
H=\< r,\om \>+h(\theta)
$$  
exhibits diffusion at all times.
\end{Main}

\bigskip

By a {\it generic set } in this paper we mean a dense $G_\delta$ set.
\begin{Main}\label{th_weak_mixing}
For any $d\geq 3$, for a generic  $\om\in \RR^d$, any
$\de>0$ and generic $h_1(\theta)$, $h_2(\theta)\in \cO^{\om}_\de$,
the flow $\Phi^t_H$ of the Hamiltonian 
$$
H(\theta,r)=\frac{1}{1+h_1(\th)}  (\< r,\om \> + h_2(\theta))
$$
is topologically weakly mixing.  
\end{Main}

Note that the energy surface in the statement above is unbounded, which drove us to treat the topological version of weak mixing rather than the classical notion.
 
 Diffusion at all times shows that the rigidity property of the
 rotator (convergence of the dynamics to identity along a subsequence
 of times) can be destroyed on all energy surfaces by a small perturbation. However, our examples are not topologically mixing, and it would be very interesting to produce examples that are topologically mixing on energy surfaces.

An interesting question concerns the possibility of similar
examples in a neighborhood of an elliptic equilibrium. 
If the components of the frequency vector at the equilibrium are all of the same sign, the energy
surfaces are bounded. We hope that the methods of \cite{FayadSaprykina15}
can be used to construct smooth examples of topologically weakly mixing  Hamiltonians 
in this context as well.

From the form of the perturbations of the rotator in Theorems A and B,
we can see that the Hamiltonians that we construct are in fact in the
closure of Hamiltonians that are conjugate to the rotator. Just like
reparametrizations of linear flows on the torus, or abelian skew
products above them, our constructions can thus be viewed {\it a
  posteriori} as particular instances of the AbC method (Approximation
by Conjugation) \cite{AK}, that is also called AK-method in reference
to Dmitry Anosov and Anatoly Katok who first introduced the method. 
The method was already used in the Hamiltonian context by Katok in \cite{K} to show the existence of integrable degenerate Hamiltonians with some particular Liouville frequencies and bounded energy surfaces that can be smoothly perturbed to become ergodic on the energy surfaces. Subsequent constructions that use the AbC method with several frequencies appeared in \cite{EFK,FayadSaprykina15,FF} to discuss the stability of elliptic equilibria and invariant tori, in particular those with Diophantine frequency vectors. The examples of \cite{K} were constructed following the usual AK-method with successive conjugations of a circle action, which gives $C^\infty$ flows that are  rigid in the sense that the dynamics converge to Identity along a subsequence of times.  In our constructions, we bypass the smoothness limitation and the rigidity  of the perturbed dynamics by resorting to the reparametrization technique of translation flows in dimension larger than $3$ used in \cite{F_mixing}. This technique exploits the Liouville phenomenon in several directions  \cite{Yoc} to avoid the Denjoy-Koksma  cancellations that appear in dimension two \cite{Katok_absence,Koch}.

\section{Notations and definitions}
To alleviate the notations, we will give the proofs for $d=3$ since there is no difference at all in the proof of the general case.

\subsection{General notations} 

- For a vector $r$, its components are denoted by $r_j$, $j=1,2,3$;  
for a vector $r_0\in \RR^3$, we denote its components by $r_{0,j}$,
$j=1,2,3$.

- For $a \in\RR$ and $l>0$, we denote
$$
I(a ,l) =a +[0,l]=[a ,a +l].
$$

- For a set $S$ in the phase space, 
let $\pr_r(S)$ stand for the orthogonal projection of $S$ onto the space
of actions ($r$-space). In the similar way introduce notations $\pr_\th(S)$, 
$\pr_{r_j}(S)$, $\pr_{\th_j}(S)$ for $j=1,\dots ,d$.

- When we write that $\frac{p}{q}$ is a rational
number or $\frac{p}{q}\in \QQ$, we assume that $q\in \NN$, $q\geq 1$,
$p\in \ZZ$, and the numbers $p$ and $q$ are relatively prime.

- Let $\om= (\om_1,\om_2,\om_3 )\in \R^3$ be a rationally independent
vector. Then, in particular, $\om_3\neq 0$, and we can rewrite $\om=
\om_3(\frac{\om_1}{\om_3}, \frac{\om_2}{\om_3},1 )$.
Without loss of generality, the constructions are performed under the assumption that 
$\om$ is normalised: $\om= (\a,\a',1 )$.

- Denote 
$
\displaystyle H_0(\theta,r)=\left< \om , r\right>=\sum_{j=1}^3 \om_jr_j;
$

- $\Phi_{H}^t(\cdot)$ stands for the flow map  at time $t$, defined by the
Hamiltonian $H$;

- Fix $c\in \RR$ and consider an energy surface $E_H(c)=\{(\theta,r)\mid
H(\theta,r) = c \}$. 
In our constructions, $\|H-H_0 \|_1$ is small, so $E_H(c)$ is uniformly 
close to $E_{H_0}(c)=\{r\mid \<r,\om \>=c\}\times \TT^3$.

\subsection{Arithmetic reminders. Yoccoz pairs of frequencies}
Denote
$$
\vertiii{ k \a }=\inf_{p\in\ZZ} |k\a +p|.
$$
For an irrational number $\a$ there exists a sequence of rational numbers
$( \frac{p_n}{q_n})_{n\geq 1} $, called the {\it convergents} of
$\a$ such that 
$$
\vertiii{ q_{n-1}\a} < \vertiii{ k\a } \quad \text{for all  }k<q_n,
$$  
and for any $n$
\begin{equation}\label{conv_est}
\frac{1}{q_n+q_{n+1}} \leq (-1)^n (q_n \a -p_n ) \leq
\frac{1}{q_{n+1}}.
\end{equation} 

\medskip

\begin{definition} Let $\om= (\a,\a',1 )\in \R^3$ where $\a$ and $\a'$ are irrational
real numbers with the corresponding sequences of convergents 
$(\frac{p_n}{q_n})_{n\geq 1}$ and
$(\frac{p_n'}{q_n'})_{n\geq 1}$. 
We say that $\om= (\a,\a',1  )\in Y$ if for all $n=0,\dots ,\infty$, 
the denominators of the convergents of $\a$ and $\a' $, respectively, satisfy: 
\begin{equation}\label{cond_Y}
e^{ q_n}\leq q_{n}'/4 , \quad
e^{ q_n'}\leq q_{n+1}/4.
\end{equation} 
\end{definition}
By \cite{Yoc}, the set $Y$ is nonempty, of cardinality continuum. This
is the set of frequencies used in Theorem \ref{th_diffusion}.

\subsection{Intervals and rectangles in an energy
  surface}  \label{S_defJ}
Here we describe the standard sets used in the construction.
In particular, {\it intervals} are defined to be 
small one dimensional curves that lie in a given energy surface 
and whose projection  onto the 5-dimensional space $(\theta,
r_1,r_2)$ is a linear segment parallel either to the $\theta_1$ axis or to the 
 $\theta_2$ axis. The coordinate $r_3$ is
defined by the requirement that the curve lies in the energy surface. 
More precisely:

- Given $s_0=(\th_0,r_0)$ and $l>0$,  
define the {\it intervals }
\begin{equation*}
\begin{aligned}
& J^{(1)}(s_0, l) = \\
& \{(\th , r)\in I(\th_{0,1},l)  \times \{\th_{0,2} \}\times
\{\th_{0,3}\}\times \{r_{0,1}\}\times \{r_{0,2}\}\times \RR \mid
H( \th,r)= H(\th_0,r_0) \},
\end{aligned}
\end{equation*}
and 
$$
\begin{aligned}
& J^{(2)}(s_0, l) = \\
& \{(\th , r)\in \{\th_{0,1} \}\times I(\th_{0,2},l) \times
\{\th_{0,3}\}\times \{r_{0,1}\}\times \{r_{0,2}\}\times \RR \mid
H(\th,r)= H(\th_0,r_0) \},
\end{aligned}
$$

- For any $s_0=(\th_0,r_0)$, $l_1>0$ and $l_2>0$, 
define the {\it rectangle} 
\begin{equation*}
\begin{aligned}
&R (s_0 , l_1, l_2) = \\ 
& \{(\th , r)\in I(\th_{0,1},l_1)  \times I(\th_{0,2},l_2) \times
\{\th_{0,3}\}\times \{r_{0,1}\}\times \{r_{0,2}\}\times \RR \mid
H(\th,r)= H(\th_0,r_0) \}.   
\end{aligned}
\end{equation*}
As before, the projection of $R(\th_0, r_0,l_1,l_2)$ onto the space
$(\th, r_1,r_2)$ is a flat rectangle parallel to $(\th_1,\th_2)$-plane; 
$r_3$ is chosen so that $R(\th_0, r_0,l_1,l_2)\subset E_H(c)$.

We say that the {\it size of the rectangle} $R (\th_0, r_0, l_1, l_2)$ is
$l_1\times l_2$.

- Given  $n$ and $s_0=(\th_0,r_0)$, a  {\it box} $B_n(s_0)\subset E_H(c)$
is defined by
\begin{equation}\label{defBn}
\begin{aligned}
B_n (s_0)=&\{ (\th ,r)\mid  \th_j\in I\left(\th_{0,j},\frac{1}{n}\right), \, j=1,2,3 ,\\ 
& r_j\in I\left(r_{0,j},\frac{1}{n}\right),\, j=1,2 ,\,  H(\th,r)= H(\th_0,r_0)\}. 
\end{aligned}
\end{equation}
These sets, having full dimension in $E_H(c)$,  will be used as test sets: 
in particular, to prove Theorem \ref{th_weak_mixing}, 
we will show that at certain times $t_n$, the image of any rectangle
$R_n \subset E_H(c)$ intersects each box $B_n\subset E_H(c)$. To do so, we will need the notion of stretching:

\medskip 

\begin{definition} Given positive $l $, $L$ and $t$, we say that 
{\it the flow map $\Phi_H^t$  is   $(1 , l , L)$-stretching}
if for any  
interval $J^{(1)}=J^{(1)}(\th_0,r_0, l)$ with $|r_0|\leq L/10$ we have:  
$$
\pr_{r_1}( \Phi^t_H(J^{(1)})) \supset [-L,L],
$$
{and the map $(\th,r) \mapsto \pr_{r_1}\left(
 \Phi^t_{H}(\th ,r)  \right)$ 
is independent of $\th_2$ where $\th=(\th_1,\th_2,\th_3)$}.
Analogously, we say that {\it the flow map
$\Phi_H^t$  
is $(2 , l , L)$-stretching } 
if for any interval $J^{(2)}=J^{(2)}(\th_0,r_0,l)$  with $|r_0|\leq
L/10 $ we have:  
$$
\pr_{r_2}(\Phi^t_H(J^{(2)})) \supset [-L,L],
$$
{and the map $(\th,r) \mapsto \pr_{r_2}\left(
 \Phi^t_{H}(\th, r)  \right)$ 
is independent of $\th_1$}.

\end{definition}
\bigskip

\section{Proofs of the main theorems}
Here we prove the main theorems modulo the technical statements, whose demonstration is deferred to the next section.

\subsection{The construction for Theorem \ref{th_diffusion}}

\bigskip

Let us fix an arbitrary vector $(\a,\a',1 )\in Y$ with 
$(\frac{p_n}{q_n})_{n\geq 1}$ and $(\frac{p_n'}{q_n'})_{n\geq 1}$
being the corresponding sequences of convergents. 
Let \begin{align}\label{eq_h}
h(\th)&= - \sum_{n=1}^\infty h_n(\th)-\sum_{n=1}^\infty h_n'(\th), \\ 
h_n(\th)&= e^{-q_n}\cos 2\pi (q_n\th_1 - p_n\th_3),\quad h_n'(\th)=
e^{-q_n'}\cos 2\pi (q_n'\th_2 - p_n'\th_3). \nonumber
\end{align}
Theorem  \ref{th_diffusion} follows from 
\begin{theorem} \label{th_diffusion1} For any $\om\in Y$, for $h$ as in \eqref{eq_h},  the Hamiltonian flow
$\Phi^t_H(\theta,r)$
defined by the Hamiltonian 
$$
H=\< r,\om \>+h(\theta)
$$  
exhibits diffusion at all times.
\end{theorem}

\bigskip

\begin{remark} This Hamiltonian can be seen as a limit of an Anosov-Katok
type construction, i.e., it has the form: 
$$
H=\lim H^{(n)},\quad H^{(n)}=H_0\circ\Psi_n\circ \dots \circ\Psi_1 ,
$$
where $\Psi_i$ are symplectic analytic coordinate changes. 
\end{remark}

\bigskip

The proof of Theorem \ref{th_diffusion1} relies on the following proposition that is proved in Section \ref{S_stretch_actions}.

\medskip

\begin{proposition} \label{prop_stretch} For any $\om\in Y$, for $h$ as in \eqref{eq_h},  the Hamiltonian flow
$\Phi^t_H(\theta,r)$
defined by the Hamiltonian 
$$
H=\< r,\om \>+h(\theta)
$$  
satisfies for each $n$:  
\begin{itemize}
\item[{\bg $(a)$}] For all $t\in [e^{q_n}, q_{n+1}/4 ]$, $\Phi^t_{H}$ is $(1, 1/q_n, q_n)$-stretching.  
  
\item[{\bg $(b)$}] For all  $t\in [ e^{q'_n}, q'_{n+1}/4]$,  
$\Phi^t_{H}$ is  $(2, 1/q'_n, q'_n)$-stretching.
\end{itemize} 
\end{proposition} 

Here we show how this proposition implies Theorem \ref{th_diffusion1}.

\bg \noindent {\it Proof of Theorem \ref{th_diffusion1}.} \bk 
Since the sequences $(q_n)$ and $(q'_n)$ satisfy (\ref{cond_Y}), 
the union of the intervals 
$\cup_{n>N} [e^{ q_n}, q_{n+1}/4] \cup [e^{ q'_n}, q'_{n+1}/4]$
contains the half-line  $t> e^{q_{N}}$. Proposition \ref{prop_stretch}
implies that for each $t \in [e^{ q_n}, q_{n+1}/4]$, $\Phi^t_{H}$ 
stretches small rectangles in the direction of $r_1$ with a large
factor, and for each $t
\in [e^{ q'_n}, q'_{n+1}/4]$, $\Phi^t_{H}$ 
stretches small rectangles in the direction of $r_2$.

Hence $\Phi^t_{H}$ exhibits stretching with an increasingly
strong factor as 
$t \to \infty$, in at least one of the two directions $r_1$ and $r_2$.
This implies the conclusion of Theorem \ref{th_diffusion1}. \hfill {\bg $\Box$}

\bigskip


\subsection{The construction for Theorem \ref{th_weak_mixing}} 
Consider  $\om=(\a,\a',1)$, and suppose that there exist sequences
$(p_n/q_n)_{n\geq 1}$ and $(p'_n/q'_n)_{n\geq 1}$ such that
\begin{equation}\label{num_est}
q_n^4 \leq q'_{n}, \quad  |q_n\a -p_n | \leq e^{-2 q'_n }, 
\quad | q'_n\a' - p'_n | \leq e^{-q'_n} .
\end{equation}
We start by observing that the set $S$ of pairs $(\a,\a')$ 
with this assumption contains a  generic set in $\RR^2$.
This implies, of course, that the set of numbers $\om =
(\om_1,\om_2,\om_3 )= 
\om_3 \, (\a,\a',1)$
such that  $(\a,\a')\in S$ and $\om_3\in \RR$, is generic in $\RR^3$.  
\medskip

\begin{Lem}\label{lem_numbers}
There exists a generic (dense $G_\de$) set $\hat S\subset S\subset \RR^2$ of pairs
$(\a,\a')$ satisfying 
the following: there exist sequences
$(p_n/q_n)_{n\geq 1}$ and
$(p'_n/q'_n)_{n\geq 1}$ of rational numbers such that estimate (\ref{num_est}) holds for
all $n$.
\end {Lem} 
 
 \medskip 
 
\bg \noindent{\it Proof. } \bk 
We want to describe the set $S$ of pairs $(\a,\a')$ such that for any $N$
there exist
$p/q$ and  $p'/q'\in\QQ$ such that  
$q > N$,  $q' \geq q^4 $, $p,p' \in \Z$,
and
$$
|\a -p/q | \leq e^{-2 q' }/q, 
\quad | \a' - p'/q' | \leq e^{-q'}/q' .
$$
The set $S$ contains the following set $\hat S$:
$$ 
\bigcap_{N\geq 1} 
\left(
\bigcup_{q\geq N} \bigcup_{p\in \Z} \bigcup_{q'\geq q^4} \bigcup_{p' \in \Z} 
\left(\frac{p}{q} - \frac{e^{-2 q }}{q},\, \frac{p}{q} + \frac{e^{-2 q } }{q}\right)\times 
\left(\frac{p'}{q'} - \frac{e^{-q'}}{q'},\, \frac{p'}{q'} + \frac{e^{-q'} }{q'}  \right)
\right),
$$
which is a countable intersection (in $N$) of open dense sets. \bg \qed \bk

\bigskip

Here we present an explicit example
of a Hamiltonian whose flow is topologically weakly mixing.
It is easy to see that such examples can be produced arbitrarily close
to $H_0$. From this, the genericity of Hamiltonians with the 
weak mixing property is obtained in the standard way.

Having fixed $\omega\in\RR^3$ as in \eqref{num_est}, let  
\begin{equation}\label{eq_phi}
\phi(\th)= 1+ \sum_{n=1}^{\infty}  q_n e^{-q'_n}  \cos 2\pi (q_n\th_1 -
p_n\th_3), 
\end{equation}
and for $\kappa_n:= q_n^2$, introduce 
\begin{align}\label{eq_tilde_h}
\tilde h(\th)&=-\sum_{n=1}^\infty  \tilde h_n(\th)
-\sum_{n=1}^\infty \tilde  h_n'(\th), \\
\tilde h_n(\th)&= \kappa_n e^{-q'_n}\cos 2\pi \kappa_n (q_n\th_1 - p_n\th_3),\quad 
\tilde h_n'(\th)=
e^{-q_n'}\cos 2\pi  (q_n'\th_2 - p_n'\th_3).\nonumber 
 \end{align} 

\medskip

Notice that since $\|\tilde h_n\|_1 \leq \kappa_n
  e^{-q'_n} e^{2\pi \kappa_n (q_n+p_n)} \leq  q_n^2 e^{4\pi q_n^3 - q'_n} $,
  assumption $q_n^4 \leq q'_n$ implies that $\| \sum_{n=1}^\infty  \tilde
  h_n(\th) \|_1<\infty$. Clearly, this implies that
  $\| \tilde h(\th)\|_1  <\infty$.

\medskip

\begin{theorem} \label{th_weak_mixing1} 
For  $\om$ as in \eqref{num_est}, $\phi$ as in \eqref{eq_phi}, 
$\tilde{h}$ as in \eqref{eq_tilde_h}, the Hamiltonian flow
$\Phi^t_{\tilde H}(\theta,r)$
defined by  
\begin{equation} \label{tH} 
\tilde H=\frac{1}{\phi(\th)}\left(\< r,\om \>+\tilde{h}(\theta)\right)
\end{equation} 
is topologically weakly mixing.  

More precisely,  for $t_n= e^{q_n'}$, $n\geq 1$, we have: 
for any two open sets $A$ and $B$ on the same energy surface
there exists $N=N(A, B)$ such that 
$$
\Phi_{\tilde H}^{t_n}(A)\cap B \neq \emptyset \quad \text{for all} \quad n\geq N.
$$
\end{theorem}

Assuming Theorem \ref{th_weak_mixing1}, we show how it yields Theorem \ref{th_weak_mixing}.

\medskip 

\bg \noindent{\it Proof of Theorem \ref{th_weak_mixing}.} \bk 

It follows from classical arguments (Cf. \cite{halmos}) that weak mixing for the flows as in Theorem  \ref{th_weak_mixing} holds for a  $G_\de$-set of functions $(h_1,h_2) \in \cO^{\om}_\de(0)^2$. It is
left to show the density of weak mixing for 
$(h_1,h_2) \in \cO^{\om}_\de(0)^2$ for a fixed $\de$.
To do this, notice that for any $\eps>0$, both $(\phi(\th) - 1)$ 
and $\tilde h(\th)$ can be
chosen $\eps$-close to zero in the fixed norm: it is enough to choose
$q_1$ large enough. Moreover,  from the proof of Theorem \ref{th_weak_mixing1} 
it follows that the same result holds true if we change 
$\phi(\th)$ and $\tilde h(\th)$ by 
$\phi(\th) + P(\th)  \in \cO^{\om}_\de(0)$ and $\tilde h(\th)+ Q(\th) \in \cO^{\om}_\de(0)$ with $P$ and $Q$ trigonometric polynomials. This implies the density of the weak mixing property.  \hfill {\bg $\Box$}

The proof of Theorem \ref{th_weak_mixing1} relies on the following two propositions that are proved in Sections \ref{S_stretch_actions} and \ref{S_stretch_angles}, respectively.

\medskip

\begin{proposition}\label{prop_stretch2} For $\om$ as in \eqref{num_est}, 
$\phi$ as in \eqref{eq_phi}, $\tilde{h}$ as in \eqref{eq_tilde_h}, 
$\tilde{H}$ as in \eqref{tH}, and  $t_n= e^{q'_n}$, we have:

\begin{itemize}
\item[{\bg $(a)$}]   $\Phi^{t_n}_{\tilde H}$ 
is $(1,\frac{1}{q_n^3},q_n )$-stretching;
\item[{\bg $(b)$}]  $\Phi^{t_n}_{\tilde H}$ is $(2,\frac{1}{q'_n} ,
q'_n )$-stretching.
\end{itemize} 
\end{proposition}

\bigskip

\begin{proposition}\label{prop_stretch_theta} For $\om$ as in \eqref{num_est}, $\phi$ as in \eqref{eq_phi}, $\tilde{h}$ as in \eqref{eq_tilde_h}, the Hamiltonian flow
$\Phi^t_{\tilde H}(\theta,r)$
defined by  
$$
\tilde H=\frac{1}{\phi(\th)}\left(\< r,\om \>+\tilde{h}(\theta)\right)
$$  
satisfies for $t_n= e^{q'_n}$ the following. For any rectangle 
$R_n:=R(\th_0, r_0, 1/q_n, 1/q_n)$ with $|r_0|\leq n$ and any box $B_n$ (see notations in Sec. \ref{S_defJ}) there exists 
a rectangle  $R'_n \subset R_n$ of 
size  $1/q_n^3 \times 1/q_n^3$ such that 
$$
\pi_\th (\Phi_{\tilde H}^{t_n}(R'_n)) \subset \pi_\th (B_n).
$$

\end{proposition}

\bigskip 

\bg \noindent {\it Proof of Theorem \ref{th_weak_mixing1}.} \bk 
Fix $R_n$ and $B_n$ as above.

By Proposition \ref{prop_stretch_theta},
there exists a rectangle  $R'_n \subset R_n$ of size $1/q_n^3 \times
1/q_n^3$ such that 
$$
\pi_\th (\Phi_{\tilde H}^{t_n}(R'_n)) \subset \pi_\th (B_n).
$$
By Proposition \ref{prop_stretch2}, we can find a rectangle $\bar R_n \subset R'_n$ such that 
$$
\pi_r (\Phi_{\tilde H}^{t_n}(\bar R_n)) \subset \pi_r (B_n).
$$
Hence 
$$
\Phi_{\tilde H}^{t_n}(\bar R_n)) \subset B_n$$
and the proof is finished. \hfill $\Box$

\section{Stretching}


\subsection{ Stretching in the action directions.}\label{S_stretch_actions}

\begin{Lem}\label{lem_model1}
Let $\frac{p}{q}$, $\frac{p'}{q'}\in\QQ$ and $\om=(\a,\a',1)$ and $a_q, b_{q'}$ satisfy 
$$
e^{-q}\geq a_{q}\geq {4|q\a-p|},\quad e^{-q'}\geq b_{q'}  \geq {4|q'\a'-p'|}.
$$  
Define
$$
H(\th ,r)= \<r,\om\> - a_q\cos 2\pi  (q\th_1 - p\th_3) -  
b_{q'}\cos 2\pi  (q'\th_2 - p'\th_3).
$$ 
Then the following holds: 
\begin{itemize}
\item[{\bg $(a)$}] For each $t\in [a_{q}^{-1}, 1/(4|q\a - p|) ]$ 
the  flow map $\Phi^t_{H}$ is $(1, 1/q, 2q)$ stretching.

\item[{\bg $(b)$}] For each $t\in [b_{q'}^{-1}, 1/(4 |q'\a'- p'|) ]$, 
the  flow map $\Phi^t_{H}$ is $(2, 1/q', 2q')$ stretching.
\end{itemize}
\end {Lem}

\medskip

\bg \noindent{\it Proof of Lemma \ref{lem_model1}.} \bk 
The Hamiltonian $H$ defines the following system of equations (we
omit $r_3$ from the considerations):
\begin{equation}\label{model_system1}
\begin{cases}
\dot \theta = \omega, \\
\dot r_1 = -2\pi q a_{q}  \sin 2 \pi (q\th_1-p\th_3) , \\
\dot r_2 = -2\pi q' b_{q'}  \sin 2 \pi (q'\th_2-p'\th_3) .\\
\end{cases}
\end{equation}
This system can be integrated explicitly: the solution with 
initial conditions $(\th(0),r(0))=(\th_0,r_0)$ satisfies
\begin{equation}
\begin{cases}
\theta = \theta_{0} +t \omega, \\
r_1 (t) = c_1 + \frac{q a_{q}}{q\a-p}
\cos 2\pi ((q\th_{0,1}- p\th_{0,3})+t (q\a - p) ), \\
r_2 (t)= c_2 + \frac{q' b_{q'}}{q'\a' - p'} \cos
2\pi ((q'\th_{0,2} - p'\th_{0,3})+t (q'\a - p') ), \\
\end{cases}
\end{equation}
where $c_k$ is a constant such that $r_{k}(0)=r_{0,k}$, $k=1,2$.
Notice that $r_1(t) = \pr_{r_1}\left( \Phi^t_{H} (r_0,\th_0) \right)$ is independent of $\th_{0,2}$, and 
$r_2(t)=\pr_{r_2} \left( \Phi^t_{H} (r,\th))\right)$ is independent of
$\th_{0,1}$.

Fix an arbitrary $t\in [a_{q}^{-1}, 1/(4(q\a - p)) ]$ and $s_0=(r_0,\th_0)$ with
$|r_0|< q/10$, and consider the interval $J^{(1)}=J^{(1)}(s_0,\frac{1}{q})$,
see notations in Section  \ref{S_defJ}.
Assume $q\a - p > 0$ (the opposite case is similar).
There exists a point $s^+=(\th^+,r_0) \in J^{(1)}$ with
$\th^+= (\th^+_{1} , \th_{0,2},\th_{0,3})$ such that
$$
q\th^+_{1} - p\th_{03}=-1/4 \ \ \text{mod}\, 1.
$$
Then $\cos2\pi(q\th^+_{1} - p\th_{03})=0$.  
Consider the trajectory of the above flow with the initial condition
$(\th(0), r(0))=(\th^+,r_0)$. For the first action component this
reads: $r_1 (0) = c_1 = r_{0,1}$, which implies
$$
|c_1|=|r_{0,1}|\leq |r_0|\leq  q/10.
$$

Assumption $t\in [a_{q}^{-1}, 1/(4(q\a - p)) ]$ implies, in particular, 
$(2\pi t (q\a - p)) \in [0, \pi/2 ]$. Using the  fact that 
$\sin(x)\geq x/2$ for all $x\in [0, \pi/2 ]$, we get
$$
\begin{aligned}
&\frac{q a_{q}}{q\a - p}\cos 2\pi ((q\th^+_{1} - p\th_{03})+t
(q\a - p) )= \\
&\frac{q a_{q}}{q\a - p} \sin 2\pi t (q\a - p) \geq \pi t q a_{q} 
\geq 3q.
\end{aligned}
$$
Since $|c_1|< q/10$, we get
$\pi_{r_1}\left( \Phi^t_{H}(s^+)\right) > 2q$.

In the same way, there is a point $s^- \in J^{(1)}$  
such that $\pi_{r_1}( \Phi^t_{H}( s^-)) <- 2q$.
The result follows by continuity. \hfill {\bg $\Box$}

\bigskip

Here we  
prove Proposition \ref{prop_stretch} that was used for the proof of Theorem  \ref{th_diffusion1}.

\medskip

\bg \noindent{\it Proof of Proposition \ref{prop_stretch}.} \bk 

Let us prove statement $(a)$, the second statement is similar. 
Since $q_n/p_n$ is a 
convergent of $\a$,  
by (\ref{conv_est}) we have for all $n$ that 
$
q_{n+1}\leq \frac{1}{|q_n\a-p_n|}
$.
Condition (\ref{cond_Y}) implies that $e^{q_n}\leq q_{n+1}/4$, so we have
$$
[\, e^{q_n}, q_{n+1}/4 \,] \subset [\,e^{q_n}, 1/(4|q_n\a-p_n|) \,]. 
$$
Consider $H_n(\th ,r)= \<r,\om\> - h_n(\th) -  h_n'(\th)$, 
where $h_n$, $h_n'$  are defined by (\ref{eq_h}). Fix $r=r_0$ with $|r_0|<q/10$.
The first component of $\Phi^t_{H_n}(\th,r_0)$, i.e., 
$r_{n,1} (\th,r_0,t):=\pr_{r_1}(\Phi^t_{H_n}(\th,r_0))$, is
given by
$$
\begin{aligned}
r_{n,1} (\th,r_0,t) = &c_{n}(\th,r_0) + \frac{q_{n}e^{-q_n} }{q_n\a- p_n}  
\cos  2\pi ( (q_n\theta_{1}-p_n\theta_{3}) + t
(q_n\a -p_n) ):=  \\
&c_n(\th,r_0) + f_n(\th,t),
\end{aligned}
$$
where $c_n(\th,r_0)$ is such that $r_{n,1} (\th,r_0,0)=r_0$. 
By Lemma  \ref{lem_model1}, for each $t\in [\,  e^{q_n}, q_{n+1}/4 \, ] $ 
the  flow map $\Phi^t_{H_n}$ is 
$(1, 1/q_n, 2q_n)$ 
stretching, i.e.,
in any interval $J^{(1)}(r_0,\th_0, 1/q_n)$ with $|r_0| \leq q_n/10$ 
there are points $s^+=(r_0,\th^+)$ and $s^-=(r_0,\th^-)$
such that
$$
r_{n,1} (s^+,t) \geq 2q_n, \quad
r_{n,1} (s^-,t) \leq -2q_n.
$$
Hence, in particular,
$$
f_n(\th^+,t)- f_n(\th^+,t) =
r_{n,1}(s^+,t)-r_{n,1}(s^+,0) \geq 2q_n- q_n/10 = 1.9 q_n.
$$

Let us show that for these $t$, $\Phi^t_{H}$ has the same stretching
properties as $\Phi^t_{H_n}$ (with $2q_n$ replaced by $q_n$). 
The first component 
of $\Phi^t_{H}(\th,r_0)$, i.e., 
$r_1 (\th,r_0,t):=\pr_{r_1}(\Phi^t_H(\th,r_0)$,  is given by the formula
$$
\begin{aligned}
r_1 (\th,r_0,t) = &c(\th,r_0) + \sum_{k=1}^{\infty} \frac{q_{k}e^{-q_k} }{q_k\a- p_k}  
 \cos  2\pi ( (q_k\theta_{1}-p_k\theta_{3}) + t
   (q_k\a - p_k) ):= \\
&c(\th,r_0) + \sum_{k=1}^{\infty} f_k(\th,t)=c(\th,r_0) +f_n(\th,t)+
C_n(\th,t) +  D_n(\th,t) ,
\end{aligned}
$$ 
where $C_n(\th,t)=\sum_{k=1}^{n-1} f_n(\th,t)$, $D_n(\th,t)=\sum_{k=n+1}^{\infty}
f_k(\th,t)$, and $c(\th,r_0)$ is a constant such that $r_1
(\th,r_0,0)=r_{0,1}$.  

Consider  $C_n(\th,t)$. By (\ref{conv_est}), $1/|q_k\a-p_k|
\leq 2q_{k+1}$ for all $k\geq 1$, so for any $\th$ we have
$$
|C_n(\th,t)-C_n(\th,0)|  \leq 2 \sum_{k=1}^{n-1}  \left| \frac{q_{k}e^{-q_k} }{q_k\a- p_k } 
   \right| 
\leq  4\sum_{k=1}^{n-1} q_{k} e^{-q_k} q_{k+1} \leq  q_{n}/100
$$ 
due to the growth condition on $q_k$.

Now consider $D_n(t)$.  By (\ref{conv_est}), for any $k\geq 1$ we
have: 
$|q_{k}\a -p_{k} |\leq \frac{1}{q_{k+1}}$. 
Hence, for any  $t\leq q_{n+1}$  we have: 
$$
\begin{aligned}
|D_n(\th,t)-& D_n(\th,0)| \leq t \, \sup |D'_n(\th,\tilde t)| \leq \\ 
& 2\pi t \sum_{k=n+1}^{\infty}  \left| \frac{q_{k}e^{-q_k} }{q_k\a- p_k } 
   \right||q_k\a - p_k | \leq 2\pi q_{n+1}   \sum_{k=n+1}^{\infty} q_{k}e^{-q_k} <1.
\end{aligned}
$$ 
This implies that for any $r_0$ with $|r_0|\leq q/10$ we have:
$$
\begin{aligned}
r_1 (s^+,t)-&r_1 (s^+,0) \geq  
 f_n(\th^+,t)- f_n(\th^+,t) - |C_n(\th,t)-C_n(\th,0)|- \\
& |D_n(\th,t)-D_n(\th,0)|
\geq 1.9 q_n- 0.01 q_n -1 >  1.5 q_n .
\end{aligned}
$$ 
Since, by assumption, 
$|r_1 (s^+,0)| =|r_0| \leq q_n/10$, we get 
$$
r_1 (s^+,t) \geq q_n.
$$
By the same argument,
$r_1 (s^-,t) \leq -q_n$. Thus,
$\Phi^t_{H}$ is $(1, 1/q_n, q_n)$ stretching. \hfill {\bg $\Box$}

\bigskip

Below we prove an analog of Proposition \ref{prop_stretch} for the
Hamiltonian $\tilde H$ of Theorem \ref{th_weak_mixing1}.
To begin with, notice that our choice of $\phi$ and $\tilde h$ implies
that the Hamiltonian system of $\tilde H$ has a particularly simple
form.

\medskip 

\begin{Lem}\label{lem_} For $\om$ as in \eqref{num_est}, $\phi$ as in \eqref{eq_phi}, $\tilde{h}$ as in \eqref{eq_tilde_h}, the Hamiltonian flow
$\Phi^t_{\tilde H}(\theta,r)$
defined by  
$$
\tilde H=\frac{1}{\phi(\th)}\left(\< r,\om \>+\tilde{h}(\theta)\right)
$$  
satisfies
$$
\begin{cases}
\dot \theta = \frac{1}{\phi(\th)} \omega, \\
\dot r_1 = \frac{-2\pi}{\phi(\th)} \sum_{n=1}^\infty q_n       
\left( C q_ne^{-q'_n} \sin 2\pi (q_n\th_1 - p_n\th_3) + 
\kappa_n^2 e^{-q'_n} \sin 2\pi \kappa_n (q_n\th_1 - p_n\th_3) \right), \\
\dot r_2 = \frac{-2\pi}{\phi(\th)}   \sum_{n=1}^\infty 
q'_n  e^{-q'_n} \sin 2\pi (q'_n\th_2 - p'_n\th_3) ,
\end{cases}
$$
where $C=\tilde H(\th ,r)$.
We omit the expression for $r_3$ since it is not used below. 
\end{Lem}

\medskip

\bg \noindent{\it Proof of Lemma \ref{lem_}.} \bk 
Recall that the value of $\tilde H(\th ,r)=\frac{1}{\phi(\th)}(\<
r,\om \> + \tilde h(\th ) ):= C$ is constant on the solutions of
the corresponding system of equations. 
For $j=1,2$ we have: $\dot r_j=- \partial_{\th_j}\tilde H$, where
$$
\begin{aligned}
\partial_{\th_j}\tilde H = -\frac{1}{\phi^2}\, \partial_{\th_j}\phi \,
(\< r,\om \> +\tilde h) 
+ \frac{1}{\phi} \,\partial_{\th_j}\tilde h 
=- \frac{1}{\phi}\, (C\, \partial_{\th_j}\phi  - \partial_{\th_j}\tilde h).
\end{aligned}
$$
Explicit substitution finishes the proof. 
\hfill {\bg $\Box$}

\bigskip 


In the next lemma (which is an analog of Lemma \ref{lem_model1}) 
we study the action components of the above
system in a simplified form: we consider only the $n$-th term in the sums
above. 
The study of the angle components is postponed
to Proposition \ref{prop_stretch_theta}.

\medskip 

\begin{Lem}\label{lem_model2}
Let $\om=(\a,\a',1)$, where $\a$ and $\a'$ are irrational.
Assume that there exist rational numbers  
$p/q$ and $p'/q'$ satisfying \eqref{num_est} with $q_n$ and $q'_n$
replaced by $q$ and $q'$, respectively.
Let $\phi(\th):\TT^3\mapsto\RR$ be a smooth function satisfying 
$
3/4 \leq |\phi(\th)|\leq 2
$
for all $\th\in \TT^3$. 
Denote $\kappa=q^2$, 
take any $C\in \RR$ with $|C|\leq q^{1/2}$, and $r_0\in\RR^3$ 
with $|r_0|\leq q/10$,
and let   $\Phi_S^t(\th_0,r_0)$ be the flow of the system 
\begin{equation}\label{eq_Sn}
\begin{cases}
\dot \theta = \frac{1}{\phi(\th)} \omega, \\
\dot r_1 = \frac{-2\pi }{\phi(\th)}   q       
\left( C q e^{-q'}  \sin 2\pi (q\th_1 - p\th_3) + 
\kappa^2 e^{-q'} \sin 2\pi \kappa (q\th_1 - p\th_3 ) \right), \\
\dot r_2 = \frac{-2\pi}{\phi(\th)}  q' e^{-q'} \sin 2\pi (q'\th_2 -
p'\th_3) 
\end{cases}
\end{equation}
with initial conditions $\Phi_S^0(\th_0,r_0)=(\th_0,r_0)$.
Then for $t= e^{q'}$ the following holds:

\begin{itemize}
\item[{\bg $(a)$}]   $\Phi^t_{S}$ 
is $(1,\frac{1}{q \kappa}, 2q )$-stretching;
\item[{\bg $(b)$}]  $\Phi^t_{S}$ is $(2,\frac{1}{q'} ,
2 q' )$-stretching.
\end{itemize} 

\end{Lem}

\medskip

\bg \noindent{\it Proof of Lemma \ref{lem_model2}.} \bk The proof is
analogous to that of Lemma \ref{lem_model1}. The only difference is that in
this case $\dot \th(t)$ is not constant. 
Let us study 
the $r_1$-component of $\Phi_S^t(\th_0,r_0)$ (the analysis of $r_2(t)$ is similar). 
Fix an arbitrary $(\th_0,r_0)$ with $|r_0|\leq q/10$ and let 
$\tilde J^{(1)}=J^{(1)} (\th_0,r_0,\frac1{\kappa q})$.   

Since $3/4 \leq |\phi(\th)|\leq 2$, the Mean Value theorem
implies that the angle variables satisfy  
$$
\th(\th_0;t)=\th_0+t \,\xi(\th_0,t) \om, 
$$
where $3/4 \leq |\xi (\th,t)|\leq 2$  for all $\th, t$. 

Given an initial condition $r_{1}(\th,r_0,0)=r_{0,1}$, the system defines
$$
\begin{aligned}
&r_1(\th, r_0;t) =    c_1 + \frac{C q^2 e^{-q'} }{q\a - p}
 \cos 2\pi (q\th_1(\th ;t) - p\th_3(\th ;t)) + \\
& \frac{q \kappa e^{-q'} }{q\a - p} \cos 2\pi \kappa  (q\th_1(\th ;t) - p\th_3(\th ;t) )
:= c_1(\th,r_0) + g_1(\th ,r_0;t) + g_2(\th ,r_0;t)
, 
\end{aligned}
$$
where $c_1(\th,r_0)$ is the constant such that
$r_{1}(\th ,r_0,0)=r_{0,1}$. 
Notice that
$g_2(\th_0,r_0,t)$ is the leading term in the above expression.   
Assume that $q\a - p >0$, the opposite case is similar.
Clearly, there exists a point $s^+=(\th^+,r_0) \in \tilde J^{(1)}$ with
$\th^+= (\th^+_{1} , \th_{02},\th_{03})$ such that
$$
\kappa (q\th^+_{1} - p \th_{03})=-1/4 \ \text{mod}\, 1.
$$ 
We will show that  $r_1(\th^+, r_0; e^{q'}) >2q$.

First consider the term $g_2$. Notice that $g_2(\th^+,r_0,0)=0$.  
For $\th(\th^+;t)=\th^+ + t \,\xi(\th^+;t) \om$ we have:
$$
\kappa( q \th_1(\th^+;t)-p \th_3(\th^+;t)) = 
\left( -1/4  +  \kappa t \xi(\th^+; t) (q\a -p) \right)  \ \text{mod}\, 1.
$$
Then for $t= e^{q'}$  we have:
$$
\begin{aligned}
g_2(\th^+,r_0; t)=&\frac{q\kappa e^{- q'}}{q\a - p} 
 \cos 2\pi \kappa ( q \th_1(\th^+;t)-p \th_3(\th^+;t)) = 
\frac{q\kappa e^{-q'}}{q\a - p} 
\sin \left( 2\pi \kappa  t \xi(\th^+, t) (q\a - p) \right) \\
&\geq \frac{q \kappa e^{-q'}}{q\a - p} \pi \kappa  t \xi(\th^+, t) (q\a - p)  
\geq q \kappa^2 e^{-q'}  t  = q^5 .
\end{aligned}
$$
We used the evident estimate 
$\sin(x)\geq x/2$ for $x\in [0, \pi/2 ]$, and
$$
2\pi \kappa  t \xi(\th^+, t) \, (q\a - p)  \leq  4 \pi \kappa 
e^{q'}  e^{-2 q'} = 4 \pi q^2 e^{-q'} \in [0, \pi/2 ].
$$

To estimate the other terms in  $r_1(\th^+,r_0;t)$, notice that 
$|c_1 + g_1(\th^+,r_0;0)| =|r_{0,1}|\leq q/10 $.
By \eqref{eq_Sn}, the derivative $\dot g_1(\th^+,r_0,t)$ satisfies:
$| \dot g_1(\th^+,r_0,t) |\leq 2\pi Cq^2 e^{-q'}\leq 2\pi q^3 e^{-q'}$.
For $t=e^{q'}$ we have: 
$\Delta g_1(\th^+,r_0,t):=|g_1(\th^+,r_0,t)- g_1(\th^+,r_0,0)| \leq
2\pi t q^2 e^{-q'} = 2\pi q^2 $, and finally:
$$
\begin{aligned}
r_1 (\th^+,&r_0; t)= c_1 + g_1(\th^+,r_0,t) + g_2(\th^+,r_0,t) = \\
&r_{0,1} +
\Delta  g_1(\th^+,r_0,t) +  g_2(\th^+,r_0,t) \geq g_2(\th^+,r_0,t)
-|r_{0,1} |-  |\Delta g_1(\th^+,r_0,t)|  \geq \\
& q^5 - q/10 - 2\pi q^3  > 2q.
\end{aligned}
$$

In the same way, there is a point $(\th^-,r_0) \in \tilde J^{(1)}$ 
such that the solution 
$r_1(\th^-,r_0; \, t)$ 
with the initial condition
$r_1(\th^-,r_0;\, 0)=r_{0,1}$
satisfies:  $r_1(\th^-,r_0;\, e^{q'}) \leq- 2q$. 
This implies that $\Phi_{S}^{t}$ is $(1,\frac{1}{\kappa q}, 2q )$ stretching. 
In the same way one verifies that $\Phi_{S}^{t}$ is $(2,\frac{1}{q'}, 2q' )$
stretching. 

\hfill {\bg $\Box$}

\bigskip


\medskip

\bg \noindent{\it Proof of Proposition \ref{prop_stretch2}.} \bk 
The proof  of Proposition \ref{prop_stretch2} follows from  Lemma \ref{lem_model2}  exactly as Proposition \ref{prop_stretch} followed from Lemma \ref{lem_model1}. 
\hfill {\bg $\Box$}

\bigskip


\subsection{Stretching in the angle directions}\label{S_stretch_angles}
In this section we prove Proposition \ref{prop_stretch_theta}.
Namely, we study the behavior of $\pi_\th(\Phi^t_{\tilde H})$, which is the
solution of the equation 
\begin{equation}\label{eq_theta}
\dot \theta = \frac{1}{\phi(\th)} \omega.
\end{equation}
Since this flow does not depend on $r$, we fix an arbitrary 
$r_0$ and omit it from the notations. To shorten the notations from
Sec. \ref{S_defJ}, we denote
$$
\Phi_\phi^t(\th_0):=\pi_\th \left(\Phi_H^t(r_0,\th_0)\right),
$$
$$
B_{\th,n}(\th_0) := \pi_\th \left(B_n(r_0,\th_0)\right),
$$
$$
R_\th(\th_0,l,l') :=\pi_\th \left( R(r_0,\th_0,l,l') \right).
$$
Notice that $R_\th(\th_0,l,l')$ is indeed a flat rectangle.

\medskip

\bg \noindent {\it Proof of Proposition \ref{prop_stretch_theta}.} \bk 
Since $\tilde{h}(\theta)$ does not depend on the $r$ variables,
the restriction of the flow of 
$\tilde H=\frac{1}{\phi(\th)}\left(\< r,\om \>+\tilde{h}(\theta)\right)
$ onto $\T^3$ is the same as that of 
$\bar{H}=\frac{1}{\phi(\th)}\< r,\om \>$. We will study the latter
one.
 
First observe that,
considering the global section $\{\theta_3=0\}$, one can see the flow of $\bar{H}$
to be equivalent to a special flow $T^t_{(\a,\a'),\varphi}$
above the translation $T_{(\a,\a')}$ on $\T^{2}$ and under a ceiling
function of the form 
\begin{equation}\label{eq_phi2}
\varphi(\th)= 1+ \sum_{n=1}^{\infty}   q_n   e^{-q'_n}  \cos 2\pi (q_n\th_1). 
\end{equation}
The phase space $M_{(\a,\a'),\varphi}$ of $T^t_{(\a,\a'),\varphi}$ is
$\T^{2} \times \R$ with the identification 
$(\theta_1,\theta_2,s+\varphi(\theta_1,\theta_2))\sim (\theta_1+\a,\theta_2+\a',s)$. 
These flows were studied in \cite{F_mixing} and the proof of the 
proposition follows from \cite{F_mixing}. For completeness, we sketch the proof here. 
Observe that, as proved in Proposition 3 of \cite{F_mixing}, for
intervals $I_n\subset \RR$ of length $(1/2-2/n)q_n^{-1}$ of the form 
$\vertiii{q_n\theta_1} \in [\frac1 n,\frac 1 2 - \frac 1 n]$ or 
$\vertiii{q_n\theta_1} \in [\frac 1 2 +\frac1 n,1 - \frac 1 n]$, 
and for $m\in [t_n/2,2t_n]$, it holds for some constant $C>0$ that for
every $\theta_1 \in I_n$ 
\begin{equation} \label{thet_stret}
C^{-1}\frac{q_n^{{2}}}{n}\leq  |\partial_{\theta_1} \varphi_m(\theta)|  \leq C{q_n^{{2}}} ,
\end{equation} 
where $\varphi_m$ stands for the $m$-th Birkhoff sum of the function $\varphi$.
The latter estimate follows from the very good rational approximation of 
$\omega_1$, since \eqref{conv_est} implies that 
$\vertiii{q_n\om_1}\leq e^{-q'_n}/4$. 
Now, the LHS of \eqref{thet_stret} implies that
$T^{t_n}_{(\a,\a'),\varphi}(I_n\times  \{\theta_2\} \times \{s\})$,
for any $\theta_2 \in \T$ and any $s\leq C$, is a union of more than
$\sqrt{q_n}$ almost vertical strips that follow the orbit of $I_n$
under the base translation $T_{(\a,\a')}$. Since $\sqrt{q_n}\gg \exp
\circ \exp (n)$, we get that  $T_{(\a,\a'),\varphi}^{t_n}(I_n\times
\{\theta_2\} \times \{s\})$ is $e^{-2n}$ dense in the space
$M_{(\a,\a'),\varphi}$. This fact, plus  the RHS of
\eqref{thet_stret}, together with the fact that there is no shear in
the $\theta_2$ direction (the ceiling function depends only on
$\theta_1$), imply that for any box $B_{\th, n}$, and for any $s\in \RR$,
there exists a rectangle  $R'_n$ of size $1/q_n^3 \times 1/q_n^3$ such that 
$$
T_{(\a,\a'),\varphi}^{t_n}(R'_n \times \{s\}) \subset B_{\th, n}.
$$ 
Going back to the original flow on $\T^d$, this implies the
requirement of the proposition. 

\hfill {\bg $\Box$}


\bigskip 




\end{document}